\documentclass[preprint,12pt]{elsarticle}

\usepackage{amsmath,amssymb,amscd}
\usepackage{amsthm}
\usepackage{subfig}
\usepackage{graphicx}
\usepackage[nodots]{numcompress}
\usepackage{verbatim}
\usepackage{fancybox}
\usepackage[usenames,dvipsnames]{color}
\usepackage{sidecap}
\usepackage{caption}
\usepackage{ifthen}
\usepackage{textcomp}
\usepackage[all]{xy}

\newcommand{\Z}{{\mathbb Z}}

\newcommand{\R}{{\mathbb R}}

\newcommand{\dem}{{\em Proof: \;}}
\newcommand{\fdem}{\hfill $\square$}

\theoremstyle{plain}
\newtheorem{teo}{Theorem}[section]
\newtheorem{lema}[teo]{Lemma}
\newtheorem{cor}[teo]{Corollary}

\theoremstyle{definition}

\theoremstyle{remark}

\begin{document}

\begin{frontmatter}

\title{Singularities in Negami's splitting formula for the Tutte polynomial}

\author{J.M.Burgos}
\address{Center for Research and Advanced Studies of the National Polytechnic Institute,\\Mathematics Department, CINVESTAV, Mexico City, Mexico, 07360.\\ \texttt{\scriptsize{Email: burgos@math.cinvestav.mx}}}

\begin{abstract}
The n-sum graph Negami's splitting formula for the Tutte polynomial is not valid in the region $(x-1)(y-1)=q$ for $q=1,2,\ldots n-1$ with the additional region $y=1$ if $n>3$. This region corresponds to (up to prefactors and change of variables) the Ising model, the $q$-state Potts model, the number of spanning forest generator and particularizations of these. We show splitting formulas for these specializations.
\end{abstract}


\end{frontmatter}


\section{Introduction}

The Tutte polynomial of a graph $G$, also known as dichromate or Tutte-Whitney polynomial, is defined as the following subgraph generating function \cite{Tutte}:
$$T(G;x,y)= \sum_{\substack{ A\subseteq G \\ V(A)= V(G) }}\ (x-1)^{\omega(A)-\omega(G)}\ (y-1)^{\omega(A)+|E(A)|-|V(G)|}$$
where $A\subseteq G$ indicates that $A$ is a subgraph of $G$ and $\omega(G)$ denotes the number of connected components of $G$.
It is the most general graph invariant that can be defined by the \textit{deletion-contraction algorithm}:
$$T(G;x,y)= T(G/ e;x,y)+T(G-e,x,y)$$
where $e$ is neither a loop (and edge with coincident endpoints) nor a bridge (an edge whose deletion increases the number of connected components), with $T(G;x,y)=x^{i}y^{j}$ if the edge set of $G$ only has $i$ bridges and $j$ loops. Here $G/ e$ and $G-e$ denote the contraction and deletion of the edge $e$ respectively. Computing the Tutte polynomial is in general an NP-hard problem \cite{Ja}.

Different specializations with respective prefactors and change of variables of the Tutte polynomial, naturally appear as classical invariants in several branches of mathematics, physics and engineering (\cite{Aigner}, \cite{Bo}, \cite{Bi},\cite{BryOx}). For example, the Jones polynomial in knot theory \cite{Jones}, the reliability polynomial in network engineering, the Ising and Potts model in statistical mechanics (\cite{Ising}, \cite{Onsager}, \cite{Potts}), the random cluster model \cite{FortuinKasteleyn}, etc. (see Table \ref{Table}).

\begin{table}
\begin{center}
    \begin{tabular}{ | l | p{7cm} |}
    \hline
    Specialization & Invariant  \\ \hline
    $xy=1$ & Jones polynomial  \\ \hline
    $y=0$ & Chromatic polynomial  \\ \hline
    $x=1,\ y\neq 1$ & Reliability polynomial  \\ \hline
    $x=0$ & Flow polynomial  \\ \hline
    $(x-1)(y-1)=2$ & Ising model  \\ \hline
    $(x-1)(y-1)=q$ & $q$-state Potts model  \\ \hline
    $y\neq 1$ & Random cluster model \\ \hline
    $y=1$ & Number of spanning forest generator \\ \hline
    $(1,1)$ & Number of spanning tree \\ \hline
    $(2,1)$ & Number of spanning forest \\ \hline
    $(1,2)$ & Number of spanning subgraph  \\ \hline
    \end{tabular}
\end{center}
\caption{Specializations of the Tutte polynomial up to prefactors and change of variables.}\label{Table}
\end{table}

\begin{figure}
\begin{center}
  \includegraphics[width=0.7\textwidth]{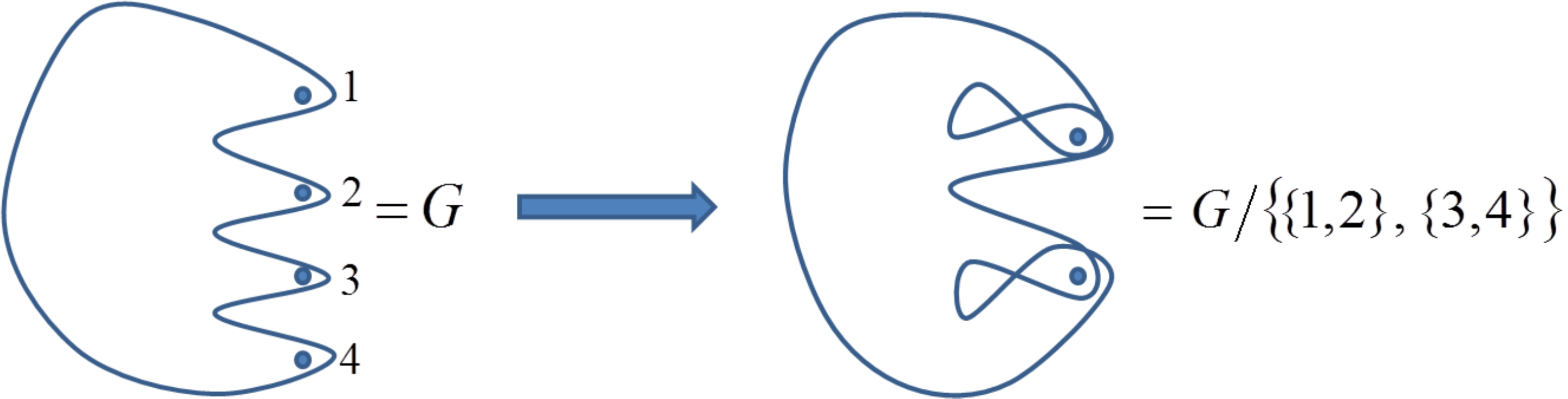}\\
  \end{center}
  \caption{Identification of vertices.}\label{Identification}
\end{figure}

Following \cite{Ne}, assume that the graph $G$ splits in subgraphs $K$ and $H$ only sharing $n$ common vertices $U=V(K)\cap V(H)$. Let $\Gamma(U)$ denote the partition lattice over $U$ and let $\mathcal{A}=\{U_{1},U_{2},\ldots U_{k}\}$ be one of these partitions. Denote by $K/\mathcal{A}$ and $H/\mathcal{A}$ the graphs obtained by identifying all vertices in each $U_{i}$ of $K$ and $H$ respectively, see Figure \ref{Identification}. The following is Negami's splitting formula for the Tutte polynomial (Corollary 4.7, iv, \cite{Ne}):
\begin{equation}\label{Tutte_Splitting}
T(G;x,y)= \sum_{\mathcal{A},\mathcal{B}\in \Gamma(U)}\ c_{\mathcal{A}\mathcal{B}}(x,y)\ T(K/\mathcal{A};x,y)\ T(H/\mathcal{B};x,y)
\end{equation}
where $c_{\mathcal{A}\mathcal{B}}(x,y)$ are rational functions of $x$ and $y$ on the field of rational numbers. A colored version of this formula was developed in \cite{Tr} and the case of Tutte polynomials of generalized parallel connections of general matroids can be found in \cite{BM}\footnote{The author is grateful to Prof.Lorenzo Traldi and Prof.Anna de Mier for these references and valuable comments.}. Explicit splitting formulas were also given in \cite{Noble} and \cite{An}. As an application of the Feferman-Vaught Theorem, the existence of splitting formulas for a wide class of graph polynomials which includes the Tutte polynomial is proved in \cite{Makowsky}. The above result is an existential theorem, it is not explicit like the others.

\begin{figure}
\begin{center}
  \includegraphics[width=1\textwidth]{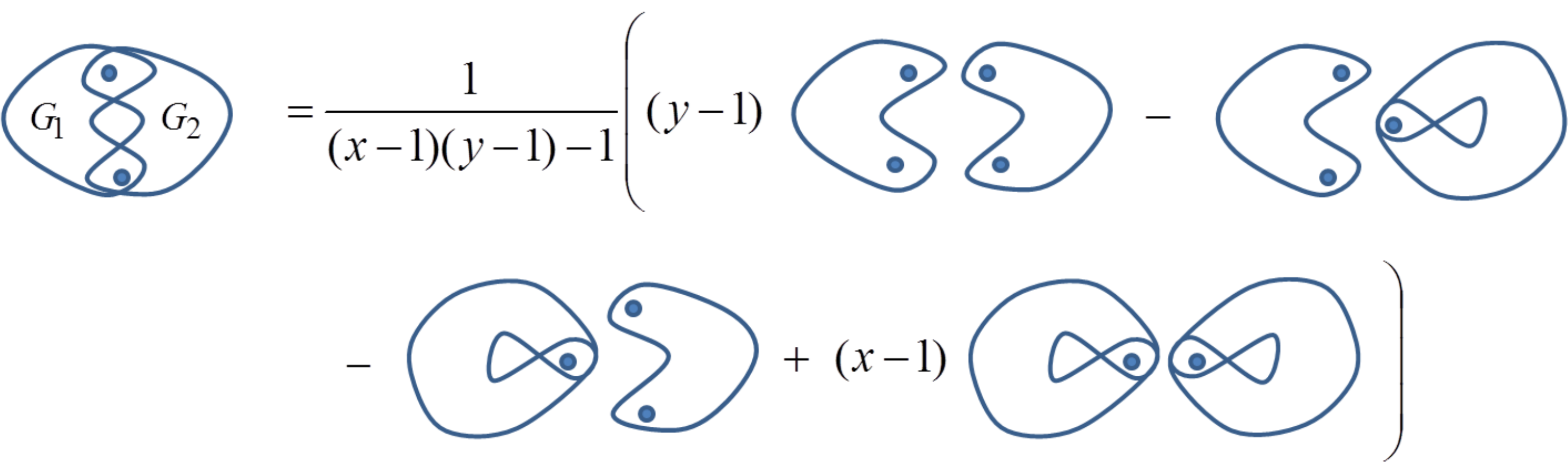}\\
  \end{center}
  \caption{Two sum Tutte polynomial splitting.}\label{nDosTutte}
\end{figure}

In view that some denominators of the coefficients $c_{\mathcal{A}\mathcal{B}}$ could annihilate restricted to certain regions, we wonder whether the formula holds for different specializations. For example, for the $2$-sum such that $H$ and $K$ are connected we have the Brylawski sum (\cite{Br}, Corollary 6.14):\footnote{As far as the author knows, this was the second known splitting formula for the Tutte polynomial after the well known factorization through an articulation point.}

\begin{eqnarray}\label{Brylawski_sum}
T(G;x,y) &=& \frac{1}{(x-1)(y-1)-1} \Big(\ (y-1)\ T(K;x,y)\ T(H;x,y) \\
\nonumber &&- T(K;x,y)\ T(H/\mathcal{A};x,y) - T(K/\mathcal{A};x,y)\ T(H;x,y) \\
\nonumber && +(x-1)\ T(K/\mathcal{A};x,y)\ T(H/\mathcal{A};x,y)\ \Big)
\end{eqnarray}
where $\mathcal{A}$ is the trivial or minimal partition of the common two vertices between $H$ and $K$. Figure \ref{nDosTutte} shows this factorization. Is clear that this formula does not hold in the region $(x-1)(y-1)=1$.

In the next section we prove the following: For the $n$-sum graph, Negami's formula \eqref{Tutte_Splitting} holds only over the region $(x-1)(y-1)\neq q$ such that $q=1,2\ldots n-1$ with the additional constraint $y\neq 1$ if $n>3$. The region where Negami's formula doesn't hold will be called the \textit{singular region}.

In this paper, in the case of the $n$-sum graph such that $H$ and $K$ are connected, we show explicit splitting formulas for the Tutte polynomial over the singular region. We also show the interesting fact that these formulas are in general not unique; i.e. The Tutte polynomial over this region splits in several different ways. For example, for $n\geq 4$, the number of spanning forest generator $T(G;x,1)$ splits according to several different formulas. These are the main results of the paper.

\section{Preliminaries}

The Negami polynomial $f(G,t,x,y)$ for a graph $G$ is defined as follows \cite{Ne}:
\begin{enumerate}
\item $f(\overline{K_{n}},t,x,y)=t^{n}$ \\
\item $f(G,t,x,y)= xf(G/e,t,x,y)+yf(G-e,t,x,y)$
\end{enumerate}
where $e$ is an edge of $G$ and $\overline{K_{n}}$ is the complement of the complete graph $K_{n}$; i.e. $n$ isolated vertices. The relationship between Negami and Tutte polynomial is the following:
\begin{equation}\label{Negami_Tutte}
f(G;(x-1)(y-1), 1, y-1)= (y-1)^{p}(x-1)^{\omega(G)}T(G;x,y)
\end{equation}

We define $\Gamma(U)$ as the set of partitions of $U$ with the following partial order: $\gamma\leq \gamma'$ if $\gamma'=\{U'_{1}, \ldots U'_{l}\}$ is a refinement of $\gamma=\{U_{1}, \ldots U_{m}\}$; i.e. For every $U'_{i}$ there is $U_{j}$ such that $U'_{i}\subseteq U_{j}$. The pair $(\Gamma(U), \leq)$ is a partition lattice. We denote by $\gamma\wedge\gamma'$ the infimum of $\gamma$ and $\gamma'$. Similarly, we denote by $\gamma\vee \gamma'$ the supremum. Consider a total order on $\Gamma(U)$ such that $\gamma_{i}\leq \gamma_{j}$ imply $i\leq j$. Define the $|\Gamma(U)|\times|\Gamma(U)|$ matrix $T_{n}$  such that its $(i,j)$-entry is $t^{|\gamma_{i}\wedge \gamma_{j}|}$ where $|\gamma|$ denotes the number of blocks of the partition.

A closer look at Negami's proof of his splitting formula\footnote{This will be done later in section \ref{y_equals_one}. It is immediate from equations \eqref{Ec_segunda_Negami} and \eqref{Ec_tercera_Negami}.} (Theorem 4.2, \cite{Ne}) shows that he actually proves the following more slightly general version:
\begin{teo}\label{Teo_Negami_Splitting}
Let $G$ be a graph obtained as a union of two graphs $K$ and $H$ sharing only the vertices $U=\{u_{1},\ldots u_{n}\}$. Let $t$ be a real number. Then,
\begin{equation}\label{Negami_Splitting}
f(G)= \sum_{\mathcal{A},\mathcal{B}\in \Gamma(U)}\ b_{\mathcal{A}\mathcal{B}}(t)\ f(K/\mathcal{A})\ f(H/\mathcal{B})
\end{equation}
such that $B_{n}(t)=\big(b_{\mathcal{A}\mathcal{B}}(t)\big)$ is a matrix verifying the relation:
\begin{equation}\label{eq_B}
T_{n}(t)B_{n}(t)T_{n}(t)=T_{n}(t)
\end{equation}
\end{teo}

If the matrix $T_{n}(t)$ is invertible for a specific value of $t$, then $B_{n}(t)= T_{n}(t)^{-1}$ and Negami's original formula is reproduced. Because of the determinant formula \cite{determinant}:
\begin{equation}\label{determinant_formula}
det\left(T_{n}(t)\right)= \prod_{\mathcal{A}\in \Gamma(U)}\ \prod_{q=0}^{|\mathcal{A}|-1}\ (t-q)
\end{equation}
the matrix $T_{n}(t)$ is non invertible for $t=0,1,\ldots n-1$ and Negami's formula does not hold. However, our generalized formulation solves this problem. For example, consider $n>1$ and the non invertible matrix $T_{n}(1)$. The matrix $B_{n}$ with one in the upper left corner and zero elsewhere is a solution of \eqref{eq_B} and gives a well defined splitting formula. Moreover, as is shown in lemma \ref{existence_sol}, for every parameter $t$ there is a solution of equation \eqref{eq_B} hence a splitting formula \eqref{Negami_Splitting} for the Negami polynomial.

Recall equation \eqref{Negami_Tutte}. Translating Negami's splitting formula \eqref{Negami_Splitting} to the the Tutte polynomial, we get the corresponding version of Corollary 4.7, iv, \cite{Ne}:
\begin{cor}\label{Splitting_non_deg}
Let $G$ be a graph obtained as a union of two graphs $K$ and $H$ sharing only the vertices $U=\{u_{1},\ldots u_{n}\}$. Then, for $(x,y)$ in the region $(x-1)(y-1)\neq 0$ we have\footnote{Recall that we are studying specializations. If we were studying the Tutte polynomial in the polynomial ring $\Z[x,y]$, the pathological region $(x-1)(y-1)=0$ wouldn't appear in the analysis.}:
$$T(G;x,y)= \sum_{\mathcal{A},\mathcal{B}\in \Gamma(U)}\ c_{\mathcal{A}\mathcal{B}}(x,y)\ T(K/\mathcal{A};x,y)\ T(H/\mathcal{B};x,y)$$
such that:
\begin{equation}\label{Tutte_coef}
c_{\mathcal{A}\mathcal{B}}(x,y)= b_{\mathcal{A}\mathcal{B}}\left((x-1)(y-1)\right)(x-1)^{\omega(K/\mathcal{A})+\omega(H/\mathcal{B})-\omega(G)}(y-1)^{|\mathcal{A}|+|\mathcal{B}|-n}
\end{equation}
where $B_{n}(t)=\big(b_{\mathcal{A}\mathcal{B}}(t)\big)$ is a solution of equation \eqref{eq_B}.
\end{cor}

See that in the original formulation by Negami, because of relation \eqref{Tutte_coef} and the determinant formula \eqref{determinant_formula}, the splitting \eqref{Tutte_Splitting} for the $n$-sum graph would hold only in the region $(x-1)(y-1)\neq q$ such that $q=0,1\ldots n-1$; i.e. Our formulation is an improvement.

\begin{figure}
\begin{center}
  \includegraphics[width=1\textwidth]{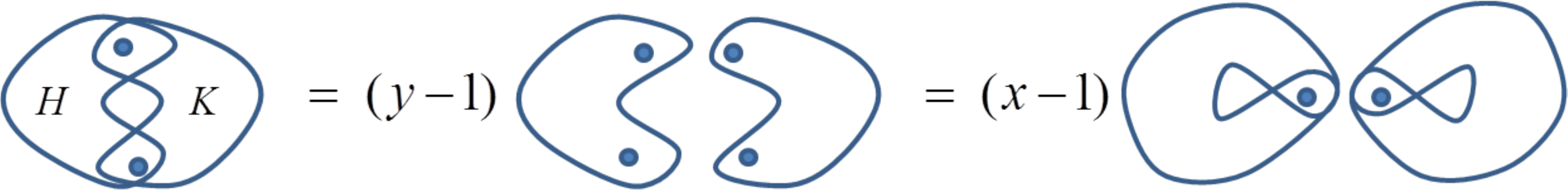}\\
  \end{center}
  \caption{Two sum factorization in the region $(x-1)(y-1)=1$.}\label{hyp1}
\end{figure}

As an example, consider the Brylawski sum \eqref{Brylawski_sum}. Along the curve $(x-1)(y-1)=1$ this splitting formula is not defined. However, Corollary \ref{Splitting_non_deg} provides the following splitting: The matrix with one in the bottom right corner and zero elsewhere is a solution of equation \eqref{eq_B} with $t=1$ hence by formula \eqref{Tutte_coef} we have the splitting:
$$T(G;x,y)=(y-1)\ T(H;x,y)\ T(K;x,y)$$
Analogously, the matrix with one in the upper left corner and zero elsewhere is a another solution and provides the splitting:
$$T(G;x,y)=(x-1)\ T(K/\mathcal{A};x,y)\ T(H/\mathcal{A};x,y)$$
where $\mathcal{A}$ is the trivial or minimal partition. These splittings hold only in the region $(x-1)(y-1)=1$ where the Brylawski sum \eqref{Brylawski_sum} is not even defined. These are illustrated in Figure \ref{hyp1}.

Recall that we are studying specializations. If we were studying the Tutte polynomial in the polynomial ring $\Z[x,y]$, the singular region wouldn't appear in the analysis for the matrix $T_{n}(t)$ is invertible in the polynomial ring\footnote{This is because the diagonal term is the higher degree term of the determinant polynomial and cannot be canceled by linear combination of the other terms.}.

What happens in the region $(x-1)(y-1)=0$? Because of equation \eqref{Negami_Tutte}, there is no relationship between Negami and Tutte polynomials when $x=1$ or $y=1$ hence we cannot assure a priori the existence of splitting formulas in this region and they have to be calculated separately. Assuming that $\omega(H)=\omega(K)=1$,  we will derive splitting formulas for this region in the following sections.

\section{Splitting formula on $x=1,\ y\neq 1$}

In the following sections we assume that $\omega(H)=\omega(K)=1$. With this condition we have that $\omega(G)=\omega(K/\mathcal{A})=\omega(H/\mathcal{B})=1$ for every pair of partitions $\mathcal{A}$ and $\mathcal{B}$ in $\Gamma(U)$ hence $\omega(K/\mathcal{A})+\omega(H/\mathcal{B})-\omega(G)=1$.

The reader may be tempted to make the following mistake: Because $\omega(K/\mathcal{A})+\omega(H/\mathcal{B})-\omega(G)=1$, there is a priori no singularity and therefore evaluating $x=1$ in \eqref{Tutte_coef} gives a trivial splitting. However, recall that the coefficients $b_{\mathcal{A}\mathcal{B}}(t)$ are rational functions and may have singularities at $x=1$.

Because the Tutte polynomial is continuous, we can just take the limit in the coefficients \eqref{Tutte_coef} and see if it defines a splitting formula; i.e. We define:
\begin{eqnarray*}
c_{\mathcal{A}\mathcal{B}}(1,y)&=& \lim_{\substack{x\to 1 \\ y\neq 1}}\ c_{\mathcal{A}\mathcal{B}}(x,y) \\
&=& \lim_{\substack{x\to 1 \\ y\neq 1}}\ (y-1)^{|\mathcal{A}|+|\mathcal{B}|-n-1}b'_{\mathcal{A}\mathcal{B}}\left((x-1)(y-1)\right)
\end{eqnarray*}
where $B'_{n}(t)=\big(b'_{\mathcal{A}\mathcal{B}}(t)\big)$ is the inverse matrix of $t^{-1}T_{n}(t)$. This is well defined for we are taking the limit and it is enough to consider a small enough reduced neighborhood of $t=1$ where the matrix is invertible. Because the inverse operation is continuous and the limit $\lim_{t\to 0} t^{-1}T_{n}(t)$ is invertible (\cite{determinant}, \cite{Bu}), we have:
\begin{equation}\label{formula_x_equals_one}
c_{\mathcal{A}\mathcal{B}}(1,y)= (y-1)^{|\mathcal{A}|+|\mathcal{B}|-n-1}b'_{\mathcal{A}\mathcal{B}}
\end{equation}
where $B'_{n}=\big(b'_{\mathcal{A}\mathcal{B}}\big)$ is the inverse matrix of the matrix $A_{n}=\big(a_{\mathcal{A}\mathcal{B}}\big)$ such that $a_{\mathcal{A}\mathcal{B}}=1$ if $|\mathcal{A}\wedge\mathcal{B}|=1$ and $a_{\mathcal{A}\mathcal{B}}=0$ otherwise. As it was expected, up to a prefactor this is the same splitting formula as the one for the Reliability polynomial studied in \cite{BR}\footnote{In that paper, the matrix $A$ was called the connectivity matrix.}:
$$R(G;p)= \sum_{\mathcal{A},\mathcal{B}\in \Gamma(U)}\ b'_{\mathcal{A}\mathcal{B}}\ R(K/\mathcal{A};p)\ R(H/\mathcal{B};p)$$
where $R(G;p)$ is the reliability polynomial of $G$.

Formula \eqref{formula_x_equals_one} shows that the singularity at $x=1$ is removable if $y\neq 1$; i.e. We have proved the following: There is an analytic continuation of the splitting coefficients $c_{\mathcal{A}\mathcal{B}}(x,y)$ to the region $x=1, y\neq 1$ such that the splitting formula \eqref{Tutte_Splitting} holds.



\section{Splitting formula on $y= 1$}\label{y_equals_one}

Direct inspection on the cases $n=1,2,3$ shows that there are no singularities at $y=1$ in the coefficients $\eqref{Tutte_coef}$ hence the splitting formula \eqref{Tutte_Splitting} holds in these cases.

However, the situation is different for $n>3$. Here, the strategy of the previous section does not work: For $n>3$ some pairs of partitions $(\mathcal{A},\mathcal{B})$ verify that $|\mathcal{A}|+|\mathcal{B}|-n-1<0$ and $b'_{\mathcal{A}\mathcal{B}}\neq 0$ hence the splitting coefficients \eqref{Tutte_coef} corresponding to these pairs when $y$ tends to one diverge:
$$\lim_{y\to 1}\ c_{\mathcal{A}\mathcal{B}}(x,y)= \lim_{y\to 1}\ (y-1)^{|\mathcal{A}|+|\mathcal{B}|-n-1}b'_{\mathcal{A}\mathcal{B}}=\infty$$
and we have no splitting formula in this region; i.e. The singularity at $y=1$ is not removable. The existence of these pairs of partitions is because otherwise the matrix $B'_{n}$ wouldn't be invertible\footnote{The proof of this fact is verbatim to the proof of the second item in Lemma \ref{existence_sol} and we will not reproduce it here.}.

We will follow Negami's strategy; i.e. Define auxiliary polynomials and derive relations between these and the contractions of the subgraphs $H$ and $K$.

In what follows, every definition or result valid for $K$ will be valid also for $H$ and the corresponding proof is verbatim. Hence we will work only with $K$. Consider the Negami polynomial expansion (Theorem 1.4 \cite{Ne}):
$$f(G;t,x,y)= \sum_{\substack{ A\subseteq G \\ V(A)= V(G) }}\ t^{\omega(A)}\ x^{|E(A)|}\ y^{|E(G)|-|E(A)|}$$
Every spanning subgraph $Y\subseteq K$ defines a partition $\mathcal{P}(Y)\in \Gamma(U)$ via the following equivalence relation: $u_{i}$ is equivalent to $u_{j}$ if they belong to the same connected component of $Y$. Define the auxiliary polynomial\footnote{This is the polynomial $A(\gamma_{i})$ in Negami's proof (Theorem 4.2 \cite{Ne}).}:
\begin{equation}\label{def_aux}
f_{\mathcal{A}}(K;t,x,y)= \sum_{\substack{ Y\subseteq K \\ V(Y)= V(K) \\ \mathcal{P}(Y)=\mathcal{A}}}\ t^{\omega(Y)-|\mathcal{A}|}\ x^{|E(Y)|}\ y^{|E(K)|-|E(Y)|}
\end{equation}

These polynomials verify\footnote{These are equations (2) and (3) in Negami's proof (Theorem 4.2 \cite{Ne}). They can also be derived by direct calculation in similar manner to \cite{BR}. These identities do not need the hypothesis $\omega(H)=\omega(K)=1$.}
\begin{equation}\label{Ec_segunda_Negami}
f(G)= \sum_{\mathcal{A}\in \Gamma(U)}\ f_{\mathcal{A}}(K)f(H/\mathcal{A})
\end{equation}

\begin{equation}\label{Ec_tercera_Negami}
f(K/\mathcal{A})=\sum_{\mathcal{B}\in \Gamma(U)}\ t^{|\mathcal{A}\wedge\mathcal{B}|}\ f_{\mathcal{B}}(K)
\end{equation}
As a consequence, we have Negami's splitting Theorem \ref{Teo_Negami_Splitting}. Define the auxiliary polynomial:
$$T_{\mathcal{A}}(K;x,y)= \sum_{\substack{ Y\subseteq K \\ V(Y)= V(G) \\ \mathcal{P}(Y)=\mathcal{A}}}\ (x-1)^{\omega(Y)-|\mathcal{A}|}\ (y-1)^{\omega(Y)+|E(Y)|-|V(K)|}$$

\begin{lema}\label{limits}
\begin{eqnarray*}
(x-1)\ T(G;x,1) &=& \lim_{\substack{ t,\zeta\to 0 \\ t/\zeta \to x-1}} f(G;t,\zeta,1)\ \zeta^{-|V(G)|} \\
T_{\mathcal{A}}(K;x,1) &=& \lim_{\substack{ t,\zeta\to 0 \\ t/\zeta \to x-1}} f_{\mathcal{A}}(K;t,\zeta,1)\ \zeta^{-|V(K)|+|\mathcal{A}|}
\end{eqnarray*}
\end{lema}
\dem
We prove only the first identity, the other is similar. Every spanning subgraph $A\subseteq G$ verifies:
$$\omega(A)+|E(A)|\geq |V(G)|$$
and the equality holds if and only if $A$ is a spanning forest with $\omega(A)$ trees. Then,
\begin{equation}\label{spanningtreeformula}
(x-1)\ T(G;x,1)= \sum_{i=1}^{|V(G)|}\ (x-1)^{i}\ S_{i}
\end{equation}
where $S_{i}$ is the number of spanning forests with $i$ trees. On the other hand we have:
\begin{eqnarray*}
f(G;t,\zeta,1) &=& \sum_{\substack{ A\subseteq G \\ V(A)= V(G) }}\ t^{\omega(A)}\ \zeta^{|E(A)|}= \sum_{i=1}^{|V(G)|}\ t^{i}\left( S_{i}\ \zeta^{|V(G)|-i}+\mathbf{O}(\zeta^{|V(G)|-i+1})\right)  \\
&=& \zeta^{|V(G)|}\ \sum_{i=1}^{|V(G)|}\ \left(\frac{t}{\zeta}\right)^{i}\left( S_{i}+\mathbf{O}(\zeta)\right)
\end{eqnarray*}
Taking the limit the result follows.
\fdem

Taking the limits of equations \eqref{Ec_segunda_Negami} and \eqref{Ec_tercera_Negami} as in lemma \ref{limits} we have\footnote{To derive equation \eqref{Ec_tercera_Tutte} from equation \eqref{Ec_tercera_Negami} we need the following fact: For every pair of partitions $\mathcal{A}, \mathcal{B}\in \Gamma(U)$ we have $|\mathcal{A}\wedge\mathcal{B}|+n -|\mathcal{A}|-|\mathcal{B}|\geq 0$. This follows from the identity:
$$|\mathcal{A}\wedge\mathcal{B}|+|\mathcal{A}\vee\mathcal{B}| \geq |\mathcal{A}|+|\mathcal{B}|$$
}:
\begin{equation}\label{Ec_segunda_Tutte}
T(G;x,1)= \sum_{\mathcal{A}\in \Gamma(U)}\ T_{\mathcal{A}}(K;x,1)\ T(H/\mathcal{A};x,1)
\end{equation}

\begin{equation}\label{Ec_tercera_Tutte}
T(K/\mathcal{A};x,1)=\sum_{\mathcal{B}\in \Gamma(U)}\ \delta^{0}_{|\mathcal{A}\wedge\mathcal{B}|+n -|\mathcal{A}|-|\mathcal{B}|}\ (x-1)^{|\mathcal{A}\wedge\mathcal{B}|-1}\ T_{\mathcal{B}}(K;x,1)
\end{equation}
where $\delta^{i}_{j}$ is the Kronecker delta. Define the matrix $L_{n}(x)=\big(l_{\mathcal{A}\mathcal{B}}(x)\big)$ whose entries are:
$$l_{\mathcal{A}\mathcal{B}}(x) =\ \delta^{0}_{|\mathcal{A}\wedge\mathcal{B}|+n -|\mathcal{A}|-|\mathcal{B}|}\ (x-1)^{|\mathcal{A}\wedge\mathcal{B}|-1}$$
For example:
$$L_{2}(x)=\left(
                         \begin{array}{cc}
                           0 & 1  \\
                           1 & x-1 \\
                         \end{array}
                       \right)$$

$$L_{3}(x)=\left(
                         \begin{array}{ccccc}
                           0 & 0 & 0 & 0 & 1	\\
                           0 & 0 & 1 & 1 & x-1 	\\
                           0 & 1 & 0 & 1 & x-1 	\\
                           0 & 1 & 1 & 0 & x-1 	\\
                           1 & x-1 & x-1 & x-1 & (x-1)^{2} 	\\
                         \end{array}
                       \right)$$

We have proved the following splitting formula:
\begin{teo}\label{Splitting2}
Let $G$ be a graph obtained as a union of two graphs $K$ and $H$ sharing only the vertices $U=\{u_{1},\ldots u_{n}\}$. Let $x$ be a real number. Then,
$$T(G;x,1)= \sum_{\mathcal{A},\mathcal{B}\in \Gamma(U)}\ d_{\mathcal{A}\mathcal{B}}(x)\ T(K/\mathcal{A};x,1)\ T(H/\mathcal{B};x,1)$$
such that $D_{n}(x)= \big(d_{\mathcal{A}\mathcal{B}}(x)\big)$ is a solution of the equation:
\begin{equation}\label{eq_D}
L_{n}(x)D_{n}(x)L_{n}(x)=L_{n}(x)
\end{equation}
\end{teo}
In lemma \ref{existence_sol} we show that there is always a solution of equation \eqref{eq_D}. For example, in the cases $n=2,3$, the unique solution $D_{n}(x)$ to equation \eqref{eq_D} is:

$$D_{2}(x)=\left(\begin{array}{cc}
1-x & 1\\
1 & 0\\
\end{array}\right)$$

$$D_{3}(x)=\frac{1}{2}\left(\begin{array}{ccccc}
(1-x)^{2} 	& 1-x & 1-x & 1-x & 2\\
1-x 		& -1 & 1 & 1 & 0\\
1-x & 1 & -1 & 1 & 0\\
1-x & 1 & 1 & -1 & 0\\
2 & 0 & 0 & 0 & 0\\
\end{array}\right)$$
In section \ref{degeneracy}, it will be shown that there is no unique solution in the case $n\geq 4$.

\subsection{Example}

Consider the case $n=4$ and the point $(x,y)=(1,1)$. Consider the following total order in the partition set $\Gamma(U)$:
$$\{\{1,2,3,4\}\}<\{\{1,2,3\},\{4\}\}<\{\{1,2,4\},\{3\}\}<\{\{1,3, 4\},\{2\}\}<\{\{2,3, 4\},\{1\}\}$$
$$<\{\{1,2\},\{3,4\}\}< \{\{1,4\},\{2,3\}\}<\{\{1,3\},\{2,4\}\}$$
$$<\{\{1,2\},\{3\},\{4\}\}<\{\{1,3\},\{2\},\{4\}\}<\{\{1,4\},\{2\},\{3\}\}$$
$$<\{\{2, 4\},\{1\},\{3\}\}<\{\{2, 3\},\{1\},\{4\}\}<\{\{2, 4\},\{1\},\{2\}\}<\{\{1\},\{2\},\{3\},\{4\}\}$$
Recall the matrix $A_{4}=\big(a_{\mathcal{A}\mathcal{B}}\big)$ such that $a_{\mathcal{A}\mathcal{B}}=1$ if $|\mathcal{A}\wedge\mathcal{B}|=1$ and $a_{\mathcal{A}\mathcal{B}}=0$ otherwise:

$$A_{4}=\left(\begin{array}{ccccccccccccccc}
1 & 1 & 1 & 1 & 1 & 1 & 1 & 1 & 1 & 1 & 1 & 1 & 1 & 1 & 1\\
1 & 0 & 1 & 1 & 1 & 1 & 1 & 1 & 0 & 0 & 1 & 1 & 0 & 1 & 0\\
1 & 1 & 0 & 1 & 1 & 1 & 1 & 1 & 0 & 1 & 0 & 0 & 1 & 1 & 0\\
1 & 1 & 1 & 0 & 1 & 1 & 1 & 1 & 1 & 0 & 0 & 1 & 1 & 0 & 0\\
1 & 1 & 1 & 1 & 0 & 1 & 1 & 1 & 1 & 1 & 1 & 0 & 0 & 0 & 0\\
1 & 1 & 1 & 1 & 1 & 0 & 1 & 1 & 0 & 1 & 1 & 1 & 1 & 0 & 0\\
1 & 1 & 1 & 1 & 1 & 1 & 0 & 1 & 1 & 1 & 0 & 1 & 0 & 1 & 0\\
1 & 1 & 1 & 1 & 1 & 1 & 1 & 0 & 1 & 0 & 1 & 0 & 1 & 1 & 0\\
1 & 0 & 0 & 1 & 1 & 0 & 1 & 1 & 0 & 0 & 0 & 0 & 0 & 0 & 0\\
1 & 0 & 1 & 0 & 1 & 1 & 1 & 0 & 0 & 0 & 0 & 0 & 0 & 0 & 0\\
1 & 1 & 0 & 0 & 1 & 1 & 0 & 1 & 0 & 0 & 0 & 0 & 0 & 0 & 0\\
1 & 1 & 0 & 1 & 0 & 1 & 1 & 0 & 0 & 0 & 0 & 0 & 0 & 0 & 0\\
1 & 0 & 1 & 1 & 0 & 1 & 0 & 1 & 0 & 0 & 0 & 0 & 0 & 0 & 0\\
1 & 1 & 1 & 0 & 0 & 0 & 1 & 1 & 0 & 0 & 0 & 0 & 0 & 0 & 0\\
1 & 0 & 0 & 0 & 0 & 0 & 0 & 0 & 0 & 0 & 0 & 0 & 0 & 0 & 0
\end{array}\right)$$
and its inverse $B'_{4}=\big(b'_{\mathcal{A}\mathcal{B}}\big)$:

$$B'_{4}=\frac{1}{6}\left(\begin{array}{ccccccccccccccc}
0 & 0 & 0 & 0 & 0 & 0 & 0 & 0 & 0 & 0 & 0 & 0 & 0 & 0 & 6\\
0 & -1 & -1 & -1 & -1 & 1 & 1 & 1 & -1 & -1 & 2 & 2 & -1 & 2 & -2\\
0 & -1 & -1 & -1 & -1 & 1 & 1 & 1 & -1 & 2 & -1 & -1 & 2 & 2 & -2\\
0 & -1 & -1 & -1 & -1 & 1 & 1 & 1 & 2 & -1 & -1 & 2 & 2 & -1 & -2\\
0 & -1 & -1 & -1 & -1 & 1 & 1 & 1 & 2 & 2 & 2 & -1 & -1 & -1 & -2\\
0 & 1 & 1 & 1 & 1 & -1 & -1 & -1 & -2 & 1 & 1 & 1 & 1 & -2 & -1\\
0 & 1 & 1 & 1 & 1 & -1 & -1 & -1 & 1 & 1 & -2 & 1 & -2 & 1 & -1\\
0 & 1 & 1 & 1 & 1 & -1 & -1 & -1 & 1 & -2 & 1 & -2 & 1 & 1 & -1\\
0 & -1 & -1 & 2 & 2 & -2 & 1 & 1 & 2 & -1 & -1 & -1 & -1 & -1 & 1\\
0 & -1 & 2 & -1 & 2 & 1 & 1 & -2 & -1 & 2 & -1 & -1 & -1 & -1 & 1\\
0 & 2 & -1 & -1 & 2 & 1 & -2 & 1 & -1 & -1 & 2 & -1 & -1 & -1 & 1\\
0 & 2 & -1 & 2 & -1 & 1 & 1 & -2 & -1 & -1 & -1 & 2 & -1 & -1 & 1\\
0 & -1 & 2 & 2 & -1 & 1 & -2 & 1 & -1 & -1 & -1 & -1 & 2 & -1 & 1\\
0 & 2 & 2 & -1 & -1 & -2 & 1 & 1 & -1 & -1 & -1 & -1 & -1 & 2 & 1\\
6 & -2 & -2 & -2 & -2 & -1 & -1 & -1 & 1 & 1 & 1 & 1 & 1 & 1 & -1
\end{array}\right)$$

\begin{figure}
\begin{center}
  \includegraphics[width=.8\textwidth]{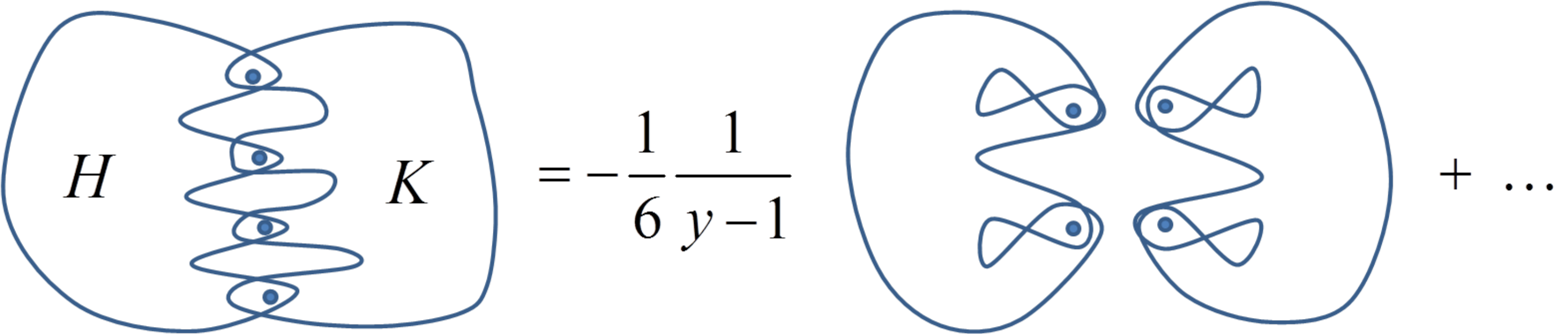}\\
  \end{center}
  \caption{Four sum Tutte polynomial splitting in the region $x=1$.}\label{Fig_Spanning_Tree}
\end{figure}

In particular, from formula \eqref{formula_x_equals_one}, the analytic continuation of formula \eqref{Tutte_coef} to the region $x=1$, we have the splitting formula:
\begin{equation}\label{Singularity_y_equals_one}
T(G;1,y)=-\frac{1}{6}\ \frac{1}{y-1}\ T(H/\{\{1,2\},\{3,4\}\};1,y)\ T(K/\{\{1,2\},\{3,4\}\};1,y)+\ldots
\end{equation}
This expression is illustrated in Figure \ref{Fig_Spanning_Tree}.

Formula \eqref{Singularity_y_equals_one} is not defined at $y=1$; i.e. Negami's formula \eqref{Tutte_Splitting} doesn't hold at the point $(x,y)=(1,1)$ in the case $n=4$. However, we still have a splitting as follows: Recall the matrix $L_{n}(x)=\big(l_{\mathcal{A}\mathcal{B}}(x)\big)$ whose entries are:
$$l_{\mathcal{A}\mathcal{B}}(x) =\ \delta^{0}_{|\mathcal{A}\wedge\mathcal{B}|+n -|\mathcal{A}|-|\mathcal{B}|}\ (x-1)^{|\mathcal{A}\wedge\mathcal{B}|-1}$$
In particular, at $x=1$ we have $l_{\mathcal{A}\mathcal{B}}(1) =\ \delta^{0}_{1+n -|\mathcal{A}|-|\mathcal{B}|}$ and the matrix is the following:

$$L_{4}(1)=\left(\begin{array}{ccccccccccccccc}
0 & 0 & 0 & 0 & 0 & 0 & 0 & 0 & 0 & 0 & 0 & 0 & 0 & 0 & 1\\
0 & 0 & 0 & 0 & 0 & 0 & 0 & 0 & 0 & 0 & 1 & 1 & 0 & 1 & 0\\
0 & 0 & 0 & 0 & 0 & 0 & 0 & 0 & 0 & 1 & 0 & 0 & 1 & 1 & 0\\
0 & 0 & 0 & 0 & 0 & 0 & 0 & 0 & 1 & 0 & 0 & 1 & 1 & 0 & 0\\
0 & 0 & 0 & 0 & 0 & 0 & 0 & 0 & 1 & 1 & 1 & 0 & 0 & 0 & 0\\
0 & 0 & 0 & 0 & 0 & 0 & 0 & 0 & 0 & 1 & 1 & 1 & 1 & 0 & 0\\
0 & 0 & 0 & 0 & 0 & 0 & 0 & 0 & 1 & 1 & 0 & 1 & 0 & 1 & 0\\
0 & 0 & 0 & 0 & 0 & 0 & 0 & 0 & 1 & 0 & 1 & 0 & 1 & 1 & 0\\
0 & 0 & 0 & 1 & 1 & 0 & 1 & 1 & 0 & 0 & 0 & 0 & 0 & 0 & 0\\
0 & 0 & 1 & 0 & 1 & 1 & 1 & 0 & 0 & 0 & 0 & 0 & 0 & 0 & 0\\
0 & 1 & 0 & 0 & 1 & 1 & 0 & 1 & 0 & 0 & 0 & 0 & 0 & 0 & 0\\
0 & 1 & 0 & 1 & 0 & 1 & 1 & 0 & 0 & 0 & 0 & 0 & 0 & 0 & 0\\
0 & 0 & 1 & 1 & 0 & 1 & 0 & 1 & 0 & 0 & 0 & 0 & 0 & 0 & 0\\
0 & 1 & 1 & 0 & 0 & 0 & 1 & 1 & 0 & 0 & 0 & 0 & 0 & 0 & 0\\
1 & 0 & 0 & 0 & 0 & 0 & 0 & 0 & 0 & 0 & 0 & 0 & 0 & 0 & 0
\end{array}\right)$$
A solution $D_4=\big(d_{\mathcal{A}\mathcal{B}}\big)$ of the equation $L_4(1)\ D_4\ L_4(1)= L_4(1)$ is the following:

$$D_4= \frac{1}{14}\left(\begin{array}{ccccccccccccccc}
0 & 0 & 0 & 0 & 0 & 0 & 0 & 0 & 0 & 0 & 0 & 0 & 0 & 0 & 14\\
0 & 0 & 0 & 0 & 0 & 0 & 0 & 0 & -3 & -3 & 4 & 4 & -3 & 4 & 0\\
0 & 0 & 0 & 0 & 0 & 0 & 0 & 0 & -3 & 4 & -3 & -3 & 4 & 4 & 0\\
0 & 0 & 0 & 0 & 0 & 0 & 0 & 0 & 4 & -3 & -3 & 4 & 4 & -3 & 0\\
0 & 0 & 0 & 0 & 0 & 0 & 0 & 0 & 4 & 4 & 4 & -3 & -3 & -3 & 0\\
0 & 0 & 0 & 0 & 0 & 0 & 0 & 0 & -4 & 3 & 3 & 3 & 3 & -4 & 0\\
0 & 0 & 0 & 0 & 0 & 0 & 0 & 0 & 3 & 3 & -4 & 3 & -4 & 3 & 0\\
0 & 0 & 0 & 0 & 0 & 0 & 0 & 0 & 3 & -4 & 3 & -4 & 3 & 3 & 0\\
0 & -3 & -3 & 4 & 4 & -4 & 3 & 3 & 0 & 0 & 0 & 0 & 0 & 0 & 0\\
0 & -3 & 4 & -3 & 4 & 3 & 3 & -4 & 0 & 0 & 0 & 0 & 0 & 0 & 0\\
0 & 4 & -3 & -3 & 4 & 3 & -4 & 3 & 0 & 0 & 0 & 0 & 0 & 0 & 0\\
0 & 4 & -3 & 4 & -3 & 3 & 3 & -4 & 0 & 0 & 0 & 0 & 0 & 0 & 0\\
0 & -3 & 4 & 4 & -3 & 3 & -4 & 3 & 0 & 0 & 0 & 0 & 0 & 0 & 0\\
0 & 4 & 4 & -3 & -3 & -4 & 3 & 3 & 0 & 0 & 0 & 0 & 0 & 0 & 0\\
14 & 0 & 0 & 0 & 0 & 0 & 0 & 0 & 0 & 0 & 0 & 0 & 0 & 0 & 0
\end{array}\right)$$
By Theorem \ref{Splitting2}, the following is a splitting formula for the spanning tree number:
$$T(G;1,1)= \sum_{\mathcal{A},\mathcal{B}\in \Gamma(U)}\ d_{\mathcal{A}\mathcal{B}}\ T(K/\mathcal{A};1,1)\ T(H/\mathcal{B};1,1)$$

\section{Degeneracy}\label{degeneracy}

In the following, $\left\lbrace
                         \begin{array}{c}
                           n 	\\
                           i 	\\
                         \end{array}
                       \right\rbrace$ denotes the Stirlng number of the second kind whose value is the number of ways to partition a set of $n$ elements into $i$ disjoint nonempty subsets. We define $\left\lbrace
                         \begin{array}{c}
                           n 	\\
                           0 	\\
                         \end{array}
                       \right\rbrace =0$ for $n\geq 1$.

\begin{lema}\label{existence_sol}
Let $t$ ($x$) be a real number. The equation \eqref{eq_B}  (\eqref{eq_D}) has solution, and the solution is unique if and only if $T_{n}(t)$ ($A(x)$) is invertible. Moreover,
\begin{enumerate}
\item The affine matrix space $\mathcal{M}_{n,t}\subseteq M_{|\Gamma(U)|}(\R)$ of solutions of the equation \eqref{eq_B} has the following dimension:
$$dim_{\R}\ \mathcal{M}_{n,t}\ = \left\{
\begin{array}{c l}
 |\Gamma(U)|^{2}-\left(\sum_{i=0}^{t}\ \left\lbrace
                         \begin{array}{c}
                           n 	\\
                           i 	\\
                         \end{array}
                       \right\rbrace \right)^{2} & \ \ \ t=0,1,2,\ldots n-1 \\
 0 & \ \ \ otherwise
\end{array}
\right.
$$
\item Equation \eqref{eq_D} has unique solution if and only if $n\leq 3$.
\end{enumerate}
\end{lema}
\dem
The matrix $T_{n}(t)$ is real and symmetric then there is an orthogonal (in particular real) matrix $O$ such that $O\ T_{n}(t)\ O^{t}$ is the diagonal matrix whose entries are the eigenvalues of $T_{n}(t)$. We can choose $O$ such that:
\begin{equation}\label{explicit_solution_matrix_fact}
O T_{n}(t)O^{t}=\left(
                         \begin{array}{cc}
                           0 & 0 	\\
                           0 & D_{k} 	\\
                         \end{array}
                       \right)
\end{equation}
where $D_{k}$ is a $k\times k$ diagonal invertible matrix. Then, all of the solutions of equation \eqref{eq_B} are the following:
\begin{equation}\label{explicit_solution_matrix_fact}
B_{n}(t)=O^{t}\left(
                         \begin{array}{cc}
                           M_{1} & M_{2} 	\\
                           M_{3} & D_{k}^{-1} 	\\
                         \end{array}
                       \right)O
\end{equation}
such that $M_{1},M_{2},M_{3}$ are arbitrary real matrices and we have the result for equation \eqref{eq_B}. See that $T_{n}(t)$ is invertible if and only if $k=|\Gamma(U)|$ and in this case, $B_{n}(t)= T_{n}(t)^{-1}$. A verbatim argument proves the result for equation \eqref{eq_D}.

\begin{enumerate}
\item Instead of an orthogonal matrix, we can choose a non orthogonal $\Lambda\in M_{|\Gamma(U)|}(\Z)$ such that\footnote{This result follows as an adaptation of the proof of the determinant formula in \cite{Bu}.}:

$$\Lambda\ T_{n}(t)\ \Lambda^{t} =$$
$$\left(
                                                   \begin{array}{ccccc}
                                                     (t-n+1)\ldots (t-1)t & 0 & \ldots & 0 & 0 \\
                                                     0 & (t-n+2)\ldots (t-1)t\ I_{\left\lbrace\substack{n \\ n-1}\right\rbrace} & \dots & 0 &0 \\
                                                     \vdots & \vdots & \ddots & \vdots & \vdots \\
                                                     0 & 0 & \ldots & (t-1)t\ I_{\left\lbrace\substack{n \\ 2}\right\rbrace} & 0 \\
                                                     0 & 0 & \ldots & 0 & t \\
                                                   \end{array}
                                                 \right)$$
and the result follows verbatim.

\item By direct calculation, the matrix $A_{n}(x)$ is invertible for $n=1,2,3$ and every $x\in\R$. It rest to show that $A_{n}(x)$ is non invertible if $n\geq 4$. Define an ordering of $\Gamma(U)$ such that $\mathcal{A}_{i}\leq \mathcal{A}_{j}$ implies $i\leq j$. With respect to this ordering, consider the upper right submatrix $A'$ of the $\mathcal{A}\mathcal{B}$ elements of $A_{n}(x)$ such that $|\mathcal{A}|=2$ and $|\mathcal{B}|=n-1$. For these entries we have:
$$|\mathcal{A}\wedge\mathcal{B}|+n -|\mathcal{A}|-|\mathcal{B}|=|\mathcal{A}\wedge\mathcal{B}|-1$$
hence all the entries not annihilated by the Kronecker delta must be one; i.e. $A'$ entries are zero or one. All of the entries to the left of $A'$ are zero for:
$$|\mathcal{A}\wedge\mathcal{B}|+n -|\mathcal{A}|-|\mathcal{B}|\geq |\mathcal{A}\wedge\mathcal{B}|\geq 1$$
and the Kronecker delta annihilates these entries. Taking the unit entry of the upper right corner of $A_{n}(x)$ as a pivot, Gauss elimination turns every entry to the right of $A'$ to zero.
Thus we have a block whose entries are zero or one of dimension $\left\lbrace
                         \begin{array}{c}
                           n 	\\
                           2 	\\
                         \end{array}
                       \right\rbrace \times \left\lbrace
                         \begin{array}{c}
                           n 	\\
                           n-1 	\\
                         \end{array}
                       \right\rbrace$. Because\footnote{These follows from the relations $\left\lbrace
                         \begin{array}{c}
                           n 	\\
                           2 	\\
                         \end{array}
                       \right\rbrace = 2^{n-1}-1$ and $\left\lbrace
                         \begin{array}{c}
                           n 	\\
                           n-1 	\\
                         \end{array}
                       \right\rbrace = \frac{n(n-1)}{2}$.} $\left\lbrace
                         \begin{array}{c}
                           n 	\\
                           2 	\\
                         \end{array}
                       \right\rbrace > \left\lbrace
                         \begin{array}{c}
                           n 	\\
                           n-1 	\\
                         \end{array}
                       \right\rbrace$ for $n\geq 4$, Gauss elimination on the block $A'$ gives a zero row and the proof is complete.
\fdem

\end{enumerate}

\section{Acknowledgments}
The author is deeply grateful to the anonymous referees for their great work, the final form of the paper is mainly due to them.

The author is grateful to \textit{Consejo Nacional de Ciencia y Tecnolog\'ia (CONACYT)}, for its \textit{C\'atedras CONACYT} program.

\end{document}